\documentstyle[amssymb]{amsart}
\def\mylabel#1{\label{#1}}
\reversemarginpar

\def\myref#1{\ref{#1}}
\reversemarginpar

\reversemarginpar

\newcommand{\rest}{{\mathord{\restriction}}}
\newcommand{\add}{\text{\normalshape\sf{add}}}
\newcommand{\cov}{\text{\normalshape\sf{cov}}}

\newcommand{\cf}{{\text{\normalshape\sf{cf}}}}

\newcommand{\dom}{{\text{\normalshape\sf {dom}}}}

\newcommand{\QED}{\hspace{0.1in} \Box \vspace{0.1in}}

\newcommand{\N}{{\cal N}}
\newcommand{\M}{{\cal M}}

\renewcommand{\>}{\rangle}

\newtheorem{theorem}{Theorem}[section]
\newtheorem{lemma}[theorem]{Lemma}

\theoremstyle{definition}
\newtheorem{definition}[theorem]{Definition}

\newcommand{\lesdot}{\mathrel{\mathord{<}\!\!\raise 
0.8 pt\hbox{$\scriptstyle\circ$}}}
\newcommand{\Proof}{{\sc Proof} \hspace{0.2in}}

\newcommand{\lft}[2]{\mathopen\ifcase#1{}\oo\or
                        \big#2\or\Big#2\else\oo\fi} 
\newcommand{\rgt}[2]{\mathclose\ifcase#1{}\oo\or
                        \big#2\or\Big#2\else\oo\fi} 

\begin{document}

\title[On the cofinality of the smallest covering ...]{On the cofinality
  of the smallest covering of 
  the real line by 
  meager sets II}
\author{Tomek Bartoszy\'{n}ski}
\address{Department of Mathematics\\
Boise State University\\
Boise, Idaho 83725, USA}
\email{{\tt tomek@@math.idbsu.edu}}
\author{Haim Judah}
\address{Abraham Fraenkel Group for Mathematical Logic\\
Department of Mathematics\\
Bar Ilan University\\
52900 Ramat Gan, Israel 
}
\email{{\tt judah@@bimacs.cs.biu.ac.il}}
\thanks{Research partially supported by 
the Israel Academy of Science, Basic Research Foundation}
\maketitle

\begin{abstract} We study the ideal of meager sets and related ideals.
  \end{abstract}
\section{Introduction}
This paper continues the line of investigation started in
\cite{BarIho89Cof} and strengthens the main result from \cite{RecOpen}. 

\begin{definition}
Suppose that ${\cal J}$ is a $\sigma$-ideal of subsets of $2^\omega$. Let 

$$\add({\cal J}) = \min\{|{\cal A}| : {\cal A} \subseteq {\cal J} \
\&\ 
\bigcup {\cal A} \not \in {\cal J}\}$$ 

and

$$\cov({\cal J}) = \min\{|{\cal A}| : {\cal A} \subseteq {\cal J} \ \&\
\bigcup {\cal A} = 2^\omega \}.$$
\end{definition}
Let $\M$ and $\N$ denote the ideals of meager and measure zero sets
respectively. 
The goal of this paper is to study the cofinality of the cardinal
$\cov(\M)$. 

Recall that for a set $H \subseteq 2^\omega \times 2^\omega $ and $x,y
\in 2^\omega$, $(H)_x = \{y \in 2^\omega : \<x,y\> \in H\}$ and 
$(H)^y = \{x\in 2^\omega : \<x,y\> \in H\} $. 

The following  definition is due to Rec{\l}aw (\cite{RecOpen}). 
\begin{definition}\label{reclaw}
  A set $X \subseteq 2^\omega$ is an $R$ set if for every Borel set $H
  \subseteq 2^\omega \times 2^\omega $, such that
  $(H)_x \in \M$ for 
  all $x \in 2^\omega$, $$\bigcup_{x \in X} (H)_x \neq 2^\omega.$$
 
\end{definition}

We will use the following representation theorem for Borel
sets with meager sections:
\begin{lemma}[Fremlin \cite{FrePa}]\mylabel{fremlin}
Suppose that $H \subseteq 2^\omega \times 2^\omega$ is a Borel set
such that $(H)_x$ is meager for all $x$. Then there
exists a sequence of Borel sets $\{G_n :n \in \omega\} \subseteq
2^\omega \times 2^\omega $ such that
\begin{enumerate}
  \item $(G_n)_x$ is a closed nowhere dense set for all $x \in
2^\omega$,
\item $H \subseteq \bigcup_{n \in \omega} G_n$.
\end{enumerate}
\end{lemma}
\Proof For completeness we present a sketch of the proof here.

Let  $\cal G$  be the family of Borel subsets  $G$  of $2^\omega \times
2^\omega$  such that $(G)_x$ is open for every $x \in 2^\omega$.

Let ${\cal J}$ be the $\sigma$-ideal of subsets of the plane generated
by Borel sets $F$ such that $(F)_x$ is closed nowhere dense for all $x$.

Consider the family $\Sigma$ of subsets  $E$  of the plane
such that there are two Borel sets  $G$, $H$  of which  $G \in\cal G$,  
$H \in {\cal J}$ and
$E\triangle G\subseteq H$.   Clearly  $\Sigma$  contains all open sets and
is closed under countable unions.  We want to show that $\Sigma$
is also closed under complements.

For a set $G \in \cal G$ let 
$$G'=\{\<x,y\> : y \text{ is an interior point of } 2^\omega \setminus
(G)_x\}.$$

Note that $(2^\omega \times 2^\omega) \setminus (G \cup G')$ is a set
whose vertical sections are  closed and nowhere dense.
It follows that in order to show  that $\Sigma $ is closed under
complements it is enough to check that $G'$ is a Borel set.

Let $\{U_n :n \in \omega\}$ be a recursive enumeration of a countable
base for the family 
of open subsets of 
$2^\omega$. 

Note that the following are equivalent:
\begin{enumerate}
\item $\<x,y\> \in G'$,
\item $ \exists n \ \lft2(y \in U_n \ \& \ \forall z \
\lft1(z \not\in U_n \ \vee\  \<x,z\> \not\in G\rgt1)\rgt2)$
(${\boldsymbol \Pi}^1_1$),
\item $\exists n \ \lft2(y \in U_n \ \& \ \forall m \ \lft1( U_n \cap U_m =
\emptyset \ \vee\  \exists z \ (z \in U_m \ \<x,z\> \not\in
G\rgt1)\rgt1)\rgt2)$ (${\boldsymbol \Sigma}^1_1$).
\end{enumerate}

That shows that $G'$ has a $\boldsymbol \Delta^1_1$ definition which
means that it is a Borel set.

So  $\Sigma$  is a $\sigma$-algebra which contains every Borel set.
Now if  $E$  is a Borel set with 
meager vertical sections then  associated  $G$  must be empty. In
particular, $E \in {\cal J}$.~$\QED$

\section{$R$ sets and  the cofinality of \normalshape\sf{cov}$({\cal M})$}

In this section we will study $R$ sets and show that $R$ sets form  an
ideal closely related to the ideal of meager sets.

\begin{definition}
  A set $X \subseteq 2^\omega$ has Rothberger's property (is a $C''$
  set) if for every sequence of open covers of $X$, $\{{\cal G}_n :n
  \in \omega\}$ there exists a sequence $\{U_n : n \in \omega\}$ with
  $U_n \in {\cal G}_n$ such that $X \subseteq \bigcup_{n \in \omega} U_n$.
\end{definition}

Rothberger's property is the topological version of strong
measure zero. We have the following:
\begin{theorem}[Fremlin, Miller \cite{FreMil88Som}]
  The following are equivalent:
  \begin{enumerate}
  \item $X \subseteq 2^\omega $ has Rothberger's property,
\item $X$ has strong measure zero with respect to every metric which
  gives $X$ the same topology.~$\QED$
  \end{enumerate}
\end{theorem}

Let $C''$ be the collection of subsets of $2^\omega$ which have
Rothberger's property. It is easy to see that $C''$ is a $\sigma$-ideal.

\begin{theorem}\mylabel{main}
  The following conditions are equivalent:
  \begin{enumerate}
  \item $X$ is an $R$ set,
  \item for every Borel function $x \leadsto f^x \in \omega^\omega$
    there exists a function $g \in 
  \omega^\omega$ such that 
$$\forall  x \in X \ \exists^\infty n \
f^x(n)=g(n).$$
\item Every Borel image of $X$ has Rothberger's property.
  \end{enumerate}
\end{theorem}
\Proof $(1) \rightarrow (2)$. Suppose that $x \leadsto f^x \in \omega^\omega$
is a Borel mapping.
Let $H=\{\<x,h\> \in X \times \omega^\omega : \forall^\infty n
\ h(n) \neq f^x(n)\}$. Clearly $H$ is a Borel set with all $(H)_x$
meager and if $g \not \in
\bigcup_{x \in X} (H)_x$ then $g$ has required properties.

\vspace{0.1in}

$(2) \rightarrow (3)$ Suppose that $Y$ is a Borel image of $X$. Let
$\{{\cal G}_n : n \in \omega\}$ be a sequence of open covers of $Y$.
Suppose that ${\cal G}_n = \{G^m_n : m \in \omega\}$ for $n \in
\omega$. For $y \in Y$ let $f^y \in \omega^\omega$ be the function
defined as
$$f^y(n) = \min\{m : y \in G^m_n\} \text{ for } n \in \omega.$$
Mapping $y \leadsto f^y$ is Borel. Thus there exists a function
$f \in \omega^\omega$ such that 
$$\forall y \in Y \ \exists^\infty n \ f^y(n)=f(n).$$
Clearly the sequence $\{G_n^{f(n)} : n \in \omega\}$ has required
properties.

\vspace{0.1in}

$(3) \rightarrow (2)$ Suppose that $x \leadsto f^x \in \omega^\omega$ is a
Borel mapping. Clearly $Y=\{f^x : x \in X\}$ has Rothberger's property in
$\omega^\omega$. Consider the families
$${\cal G}_n = \{G^m_n : m \in \omega\} \text{ where } G^m_n = \{h \in
\omega^\omega : h(n)=m\} .$$
The selector chosen for this sequence of coverings immediately gives us
the function which meets every $f^x$ in at least one point.

To get the function we are looking for we have to 
split $\omega$ into infinitely many pieces and apply the
above construction to each one of them.
 
\vspace{0.1in}

$(2) \rightarrow (1)$  
We will need several lemmas. To avoid repetitions let us define:
\begin{definition}
  Suppose that $X \subseteq 2^\omega$. $X$ is nice if for every  Borel function
  $x \leadsto f^x \in \omega^\omega$ there exists a function $g \in
  \omega^\omega$ such that 
$$\forall  x \in X \ \exists^\infty n \
f^x(n)=g(n).$$
\end{definition}
We will show first that:
\begin{lemma}\label{lem1}
  Suppose that $X$ is nice.
Then for every Borel function 
$x \leadsto \<Y^x,f^x\> \in [\omega]^\omega \times \omega^\omega$
there exists $g \in  
\omega^\omega$ such that  
$$\forall x \in X \ \exists^\infty n \in Y^x \ f^x(n)=g(n).$$
\end{lemma}
\Proof
Suppose that a Borel mapping $x \leadsto \<Y^x,f^x\>$ is given.
Let $y^x_n$ denote the
$n$-th element 
of $Y^x$ for $x \in X$.
For every $x \in X$ define a function $h^x$ as
follows:
$$h^x(n) = f^x\rest
\left\{y^x_0,\ y^x_1, \ldots,
y^x_n\right\} \hbox{ for } n \in \omega .$$
Since the mapping $x \leadsto h^x$ is Borel and functions $h^x$ can be
coded as elements of 
$\omega^{\omega}$ there is a function $h$ such that
$$\forall x \in X \ \exists^{\infty}n \ h^x(n) =
h(n) .$$
Without loss of generality we can assume that $h(n)$ is a function
from an $n+1$-element subset of $\omega$ into $\omega$. 

Define $g \in \omega^{\omega}$ in the following way.
Choose inductively 
$$z_{n} \in  \dom\lft1(h(n)\rgt1)\setminus \left\{z_{0},
z_{1},\ldots, 
z_{n-1}\right\} \text{ for } n \in \omega.$$
Let $g$ be any function such that
$g(z_{n}) = h(n)(z_{n})$ for $n \in \omega$.

We show that the function $g$ has the required properties.
Suppose that $x \in X$. Notice that the equality $h^x(n)= h(n)$ 
implies that
$$f^x(z_{n}) = g(z_{n}) \hbox{ and } z_{n} \in Y^x .$$
That finishes the proof since $h^x(n) = h(n)$ for infinitely
many $n \in \omega$.~$\QED$

\begin{lemma}\label{lem1a}
  Suppose that $X$ is nice. Then for every Borel mapping $x \leadsto
  f^x \in \omega^\omega$ there exists an increasing sequence $\{n_k :
  k \in \omega\}$ such that 
$$\forall x \in X \ \exists^\infty k \ f^x(n_k)<n_{k+1}.$$
\end{lemma}
\Proof Suppose that the lemma is not true and
let $x \leadsto f^x$ be the witness. Without loss of generality we can
assume that $f^x$ is  increasing for all $x \in X$. To get a
contradiction we will define a Borel mapping $x \leadsto g^x \in
\omega^\omega$ such that $\{g^x: x \in X\}$ is a dominating family.
That will contradict the assumption that $X$ is nice.

Define
$$g^x=\max\{\underbrace{f^x\circ f^x\circ \dots \circ f^x}_{j \text{
    times }}(i) : i,j \leq
  n\} \text{ for } n \in \omega.$$
Suppose that $g \in \omega^\omega$ is an increasing function. By the
assumption there exist $x \in X$ and $k_0$ such that 
$$\forall k\geq k_0 \ f^x\lft1(g(k)\rgt1) \geq g(k+1).$$
In particular,
$$\forall k \geq g(k_0) \ g(k) \leq g^x(k)$$
which finishes the proof.~$\QED$

\begin{lemma}\label{lem2}
  Suppose that $X$ is nice. Then for every Borel mapping $x \leadsto
  Y^x \in [\omega]^\omega$ there exists a set $Y=\{u_n: n \in
  \omega\}$ such that 
  \begin{enumerate}
    \item $u_{n+1} \geq u_n+2$ for all $n$,
    \item $\forall x \in X \ |Y \cap Y^x |=\boldsymbol\aleph_0$.
  \end{enumerate}
\end{lemma}
\Proof
By applying \myref{lem2}, we can find an increasing function $f \in
\omega^\omega$ such that 
$$\forall x \in X \ |Y^x \setminus \{f(n): n \in \omega\}| =
\boldsymbol\aleph_0.$$
Let $A_0 = \{2k : k \in \omega\}$ and $A_1 = \{2k+1 : k \in \omega\}$.
Define Borel mapping $x \leadsto \<Z^x,g^x\> \in [\omega]^\omega
\times 2^\omega$ as follows: $Z^x = \dom(g^x)$ and for $n \in \omega$,
\begin{equation*}g^x(n)=\begin{cases}
 0 & \text{if } Y^x \cap \lft1(f(n),f(n+1)\rgt1) \cap A_0 \neq
  \emptyset,\\
1  & \text{if } Y^x \cap \lft1(f(n),f(n+1)\rgt1) \cap A_1 \neq
  \emptyset\\
& \text{undefined otherwise.} 
\end{cases}
\end{equation*}
Note that the first two conditions of this definition are not
exclusive. We use either value when that happens.

By \ref{lem2}, there exists a function $h \in 2^\omega$ such that 
$$\forall x \in X \ \exists^\infty n \in Z^x \ h(n)=g^x(n).$$
Define
$$Y = \bigcup_{n \in \omega} \lft1((f(n),f(n+1)\rgt1) \cap A_{h(n)}.$$
It is clear that $Y$ has required properties.~$\QED$
 
We now return to  the proof that $(2)$ implies $(1)$.
Let $H \subseteq 2^\omega\times 2^\omega$ be a Borel set with all
fibers $(H)_x$ meager. Using \myref{fremlin}, find Borel sets $\{F^n :
n \in \omega\}$ such 
that 
$H \subseteq \bigcup_{n \in \omega} F^n$ and $(F^n)_x$ is closed
nowhere dense for $x \in 2^\omega$, $n \in \omega$.

For $x \in X$, define
$$s^{x}_{n} = \min \left\{s \in 2^{<\omega}: \forall t \in 2^{<n} \
\forall j \leq n
\ [t^{\frown}s] \cap (F^j)_x = \emptyset \right\} \hbox{ for }  n \in
\omega  .$$
(where the minimum is taken with respect to some enumeration of
$2^{<\omega}$).

By \ref{lem1a},
  there exists a sequence $\{n_k: k \in \omega\}$ such that
  \begin{enumerate}
    \item $n_{k+1} > \sum_{i=0}^k n_i$, for all $k$,
    \item $\forall x \in X \ \exists^\infty n \ |s^x_{n_k}|<n_{k+1}$.
  \end{enumerate}

For $x \in X$ let $Z^x = \{k : |s^x_{n_k}|<n_{k+1}\}$. By \ref{lem1},
there exists a sequence 
$\<s_{k} : k \in \omega \>$ 
such that 
$$\forall x \in X \ \exists^{\infty} k \in Z^x \ s^{x}_{n_k} =
s_{k}. $$
Without loss of generality we can assume that $|s_k|<n_{k+1}$ for all
$k$.
Define mapping $x \leadsto Y^x$ by
$$Y^x = \{k\in Z^x : s_k=s^x_{n_k}\}.$$
Let $Y$ be a set obtained by applying \ref{lem2} to this family.
Define
$$z = s_{l_0}\!^\frown s_{l_1}\!^\frown s_{l_2}\!^\frown \dots ,
\text{ where } l_0<l_1<l_2 \dots \text{ is the increasing enumeration of
  $Y$.}$$
Note that if $l_{k+1} \in Y \cap Y^x$ then
$$|s_{l_0}\!^\frown s_{l_1}\!^\frown \dots^\frown s_{l_k}| < \sum_{j
  \leq l_k} n_{l_j+1} <n_{l_k+2} \leq n_{l_{k+1}}$$
and $$s_{l_{k+1}}=s^x_{n_{l_{k+1}}}.$$
In other words, $[s_{l_0}\!^\frown s_{l_1}\!^\frown \dots^\frown
s_{l_{k+1}}] \cap (F^{i})_x = \emptyset$ for $i \leq l_{k+1}$.
But that clearly implies that $z \not \in \bigcup_{n \in \omega}
(F^n)_x$, which finishes the proof.~$\QED$

As a corollary we get the following:
\begin{theorem}[Rec{\l}aw \cite{RecOpen}]
Every Luzin set is an R set.
\end{theorem}
\Proof It is well known that every Borel image of a Luzin set has
Rothberger's property.~$\QED$

\begin{theorem}\mylabel{ineq}
  $\cf\lft1(\cov(\M)\rgt1) \geq \add(R) \geq \add(C'') \geq \add(\N).$ 
\end{theorem}
\Proof
We will start with the following lemma:
\begin{lemma}
There exists a Borel set $H \subseteq 2^\omega\times2^\omega $ with
such that all $(H)_x$ are meager and 
for every set $A \in \M$ there exists $x \in 2^\omega$ such that $A
\subseteq (H)_x$.   
\end{lemma}
\Proof For $z \in 2^\omega$ and $f \in \omega^\omega$ define 
a set 
$$B(z,f)=\{t \in 2^\omega : \forall^\infty n\ \exists j \in
\lft1[\widetilde{f}(n), \widetilde{f}(n+1)\rgt1) \ t(j) \neq z(j)\},
\text{ where }\widetilde{f}(n) = \sum_{k=0}^n f(k).$$
It is well known (see \cite{BJbook}) that the family $\{B(z,f) : z \in
2^\omega, f \in \omega^\omega\}$ is a basis of $\M$. 

Fix a Borel surjection $t \leadsto \<z_t, f_t\> \in 2^\omega \times
\omega^\omega $ and let $H$ be a set such that 
$(H)_t = B(z_t,f_t)$ for all $t \in 2^\omega$.~$\QED$

Suppose that ${\cal A} \subseteq \M$ is a family of meager
sets of size $\cov(\M)$ which covers $2^\omega$.
For $F \in {\cal A}$ let $x_F \in 2^\omega $ be such that $F \subseteq
(H)_{x_F}$. Suppose that $\cf(\cov(\M))=\kappa$. It follows that $X =
\{x_F : F \in {\cal A}\}$ is the union of $\kappa$ many sets
$X_\alpha$ of size smaller than $\cov(\M)$. Clearly each set $X_\alpha
\in R$ and $X \not\in R$. Thus $\kappa \geq \add(R)$ which proves the first
inequality. 

\vspace{0.1in}

The second inequality follows immediately from \myref{main}.

\vspace{0.1in}

The third inequality is due to Carlson. We will prove it here for
completeness. 
We will use the following fact (see \cite{BJbook} or \cite{Bar84Add}):
\begin{theorem}\mylabel{add}
  $\add(\N)$ is the smallest size of the family $F \subseteq
  \omega^\omega$ such that there is no function $S : \omega
\longrightarrow [\omega]^{<\omega}$ with $|S(n)| \leq n$ for all $n$,
such that 
$$\forall f \in F \ \forall^\infty n \ f(n) \in S(n).~\QED$$
\end{theorem}

Suppose that $\{X_\alpha : \alpha < \kappa < \add(\N)\}$ is a family
of $C''$ sets. We will show that $X = \bigcup_{\alpha < \kappa}
X_\alpha $ is a $C''$ set.

Let $\{{\cal G}_n : n \in \omega\}$ be a sequence of open coverings of
$X$.
Assume that ${\cal G}_n = \{U^n_m : m \in \omega\}$ for $n \in
\omega$.
Let $r(n) = 1+2+ \cdots +(n-1)$. Define for $n \in \omega$ the family
$$\widetilde{{\cal G}}_n = \left\{\widetilde{U}^n_s : s \in
\omega^{[r(n),r(n+1))}\right\}, \text{ where } \widetilde{U}^n_s =
\bigcap_{j=r(n)}^{r(n+1)} U^j_{s(j)} . $$
Note that $\widetilde{\cal G}_n$ is also a cover of $X$.
By the assumption, for every $\alpha < \kappa$ there exists a function
$f_\alpha$ such that 
$\{\widetilde{U}^n_{f_\alpha(n)} : n \in \omega\} $ is a covering of
$X_\alpha$. 

Since $\kappa<\add(\N)$, by \myref{add}, there exists a function $S: \omega
\longrightarrow \omega^{<\omega}$ such that 
$$\forall \alpha < \kappa\ \forall^\infty n \ f_\alpha(n) \in S(n).$$
Without loss of generality we can assume that $S(n)$ consists of at
most $n$
sequences of length $n$. Let $f$ be a function which agrees at least
once on the $n$-element interval $\lft1[r(n),r(n+1)\rgt1)$ with each
of the $n$ functions in $S(n)$. 
Clearly $\{U^n_{f(n)} : n \in \omega\}$ is the
covering of $X$ we were looking for.~$\QED$

{\bf Remarks and questions}

(1) Is $\cf\lft1(\cov(\M)\rgt1) \geq \add(\M)$? In other words, is
$\cf\lft1(\cov(\M)\rgt1) \geq {\frak b}$ ?

(2) The first inequality in \myref{ineq} also holds when we replace
category by measure. Unfortunately we do not know if $R$ sets defined
for measure form an ideal. (For more see \cite{Paw92Sie} and
\cite{BarJud93Add}).  

(3) It is consistent that $\add(\M) > \add(C'')$
(\cite{JudShel89Sig}) 

(4) Is it consistent that $\add(C'') > \add(\N)$? $\add(R) > \add(\N)$?


\ifx\undefined\bysame
\newcommand{\bysame}{\leavevmode\hbox to3em{\hrulefill}\,}
\fi

\end{document}